\documentclass[11pt]{amsart2000}
\usepackage{amssymb}
\newtheorem{theorem}{Theorem}[section]

\newtheorem{lemma}[theorem]{Lemma}

\theoremstyle{definition}

\theoremstyle{remark}

\numberwithin{equation}{section}
\newcommand\CH{$\mathsf{CH}$}
\newcommand\from{\colon}
\newcommand\cont{\mathfrak c}
\newcommand\con{\subseteq}
\def\CL#1{\overline{#1}}
\def\diam{\operatorname{diam}}
\def\interval{{\mathbb I}}
\def\invlim{\varprojlim}
\def\orpr#1#2{\langle #1, #2 \rangle}
\def\leukfrac#1#2{\leavevmode
\kern.1em
\raise.9ex\hbox{\the\scriptfont0 ${}_#1$}
\hskip -1pt
\kern-.1em
/\kern-.15em
\lower.10ex\hbox{\the\scriptfont0 ${}_#2$}}

\begin{document}
\setlength{\unitlength}{0.01in}
\linethickness{0.01in}
\begin{center}
\begin{picture}(474,66)(0,0) 
\multiput(0,66)(1,0){40}{\line(0,-1){24}}
\multiput(43,65)(1,-1){24}{\line(0,-1){40}}
\multiput(1,39)(1,-1){40}{\line(1,0){24}}
\multiput(70,2)(1,1){24}{\line(0,1){40}}
\multiput(72,0)(1,1){24}{\line(1,0){40}}
\multiput(97,66)(1,0){40}{\line(0,-1){40}} 
\put(143,66){\makebox(0,0)[tl]{\footnotesize Proceedings of the Ninth Prague Topological Symposium}}
\put(143,50){\makebox(0,0)[tl]{\footnotesize Contributed papers from the symposium held in}}
\put(143,34){\makebox(0,0)[tl]{\footnotesize Prague, Czech Republic, August 19--25, 2001}}
\end{picture}
\end{center}
\vspace{0.25in}
\setcounter{page}{347}
\author{Jan van Mill}
\title[Loc. conn. continuum w/o convergent sequences]{A locally connected
continuum without convergent sequences}
\address{Faculty of Sciences\\
Division of Mathematics and Computer Science\\
Vrije Universiteit\\
De Boelelaan 1081\\
1081 HV Amsterdam\\
The Netherlands}
\email{vanmill@cs.vu.nl}
\thanks{Reprinted from
Topology and its Applications,
in press,
Jan van Mill,
A locally connected continuum without convergent sequences,
Copyright (2002),
with permission from Elsevier Science \cite{vm}.}
\thanks{Jan van Mill,
{\em A locally connected continuum without convergent sequences},
Proceedings of the Ninth Prague Topological Symposium, (Prague, 2001),
pp.~347--352, Topology Atlas, Toronto, 2002}
\begin{abstract}
We answer a question of Juh\'asz by constructing under \CH\ an example
of a locally connected continuum without nontrivial convergent sequences.
\end{abstract}
\maketitle

\section{Introduction}

During the Ninth Prague Topological Symposium, Juh\'asz asked whe\-ther
there is a locally connected continuum without nontrivial convergent
sequences. 
This question arose naturally in his investigation in~\cite{juhasz2001}
with Gerlits, Soukup, and Szentmikl\'ossy on characterizing continuity in
terms of the preservation of compactness and connectedness. 
The aim of this note is to answer this question in the affirmative under
the Continuum Hypothesis (abbreviated: \CH).

Fedorchuk~\cite{vitaly77} constructed a consistent example of a compact
space of cardinality $\cont$ containing no nontrivial convergent
sequences. See also van Douwen and Fleissner~\cite{metfleissner} for a
somewhat simpler construction under the Definable Forcing Axiom. These
constructions yield zero-dimensional spaces. As a consequence, our
construction has to be somewhat different. As in~\cite{vitaly77}
and~\cite{metfleissner}, we `kill' all possible nontrivial convergent
sequences in a transfinite process of length $\omega_1$. However, our
`killing' is done in the Hilbert cube $Q=\prod_{n=1}^\infty [-1,1]_n$
instead of the Cantor set.

For all undefined notions, see \cite{engelking:gentop} and \cite{vm:book}.

\section{The Hilbert cube}\label{Hilbert}

A \emph{Hilbert cube} is a space homeomorphic to $Q$. 
Let $M^Q$ denote an arbitrary Hilbert cube.

A closed subset $A$ of $M^Q$ is a \emph{$Z$-set} if for every 
$\varepsilon>0$ there is a continuous function 
$f\from M^Q\to M^Q\setminus A$ which moves the points less than 
$\varepsilon$. 
It is clear that a closed subset of a $Z$-set is a $Z$-set.
We list some other important properties of $Z$-sets.
\begin{enumerate}
\item[(1)] Every singleton subset of $M^Q$ is a $Z$-set.
\item[(2)] A countable union of $Z$-sets is a $Z$-set provided it is closed.
\item[(3)] A homeomorphism between $Z$-sets can be extended to a
homeomorphism of $M^Q$.
\item[(4)] If $X$ is compact and $f\from X\to M^Q$ is continuous then $f$ can
be approximated arbitrarily closely by an imbedding whose range is a
$Z$-set.
\end{enumerate}
See~\cite[Chapter~6]{vm:book} for details.

Observe that by (1) and (2), every nontrivial convergent sequence with its
limit is a $Z$-set in $M^Q$.

A \emph{near homeomorphism} between compacta $X$ and $Y$ is a continuous
surjection $f\from X\to Y$ which can be approximated arbitrarily closely
by homemorphisms. This means that for every $\varepsilon > 0$ there is a
homeomorphism $g\from X\to Y$ such that for every $x\in X$ we have that
the distance between $f(x)$ and $g(x)$ is less than $\varepsilon$.

A closed subset $A\con M^Q$ has \emph{trivial shape} if it is contractible
in any of its neighborhoods. 
A continuous surjection $f$ between Hilbert cubes $M^Q$ and $N^Q$ is
\emph{cell-like} provided that $f^{-1}(q)$ has trivial shape for every
$q\in N^Q$. 
The following fundamental result is due to Chapman~\cite{chapman:lectures}
(see also~\cite[Theorem 7.5.7]{vm:book}).

\begin{enumerate}
\item[(5)] Let $f\from M^Q\to N^Q$ be cell-like, where $M^Q$ and $N^Q$
are Hilbert cubes.
Then $f$ is a near homeomorphism.
\end{enumerate}

It is easy to see that if $f\from M^Q\to N^Q$ is a near homeomorphism 
between Hilbert cubes then $f$ is cell-like. 
So within the framework of Hilbert cubes the notions `near homeomorphism' 
and `cell-like' are equivalent.

A continuous surjection $f$ between Hilbert cubes $M^Q$ and $N^Q$ is 
called a $Z^*$-map provided that for every $Z$-set $A\con N^Q$ we have
that $f^{-1}[A]$ is a $Z$-set in $M^Q$.

\begin{lemma}\label{eerstelemma}
Let $M^Q$ and $N^Q$ be Hilbert cubes, and let $f\from M^Q\to N^Q$ be a
continuous surjection for which there is a $Z$-set $A\con M^Q$ which
contains all nondegenerate fibers of $f$. 
Then $f$ is a $Z^*$-map.
\end{lemma}

\begin{proof}
Let $B\con N^Q$ be an arbitrary $Z$-set, and put $B_0= B\setminus f[A]$.
Write $B_0$ as $\bigcup_{n=1}^\infty E_n$, where each $E_n$ is compact.
It follows from~\cite[Theorem 7.2.5]{vm:book} that for every $n$ the
set $f^{-1}[E_n]$ is a $Z$-set in $M^Q$. As a consequence,
$$
f^{-1}[B] \con A\cup\bigcup_{n=1}^\infty f^{-1}[E_n]
$$
is a countable union of $Z$-sets and hence a $Z$-set by (2).
\end{proof}

\begin{theorem}\label{eerstestelling}
Let $(Q_n,f_n)_n$ be an inverse sequence of Hilbert cubes such that every
$f_n$ is cell-like as well as a $Z^*$-map. 
Then
\begin{enumerate}
\item[(A)] $\invlim (Q_n,f_n)_n$ is a Hilbert cube.
\item[(B)] The projection $f^{\infty}_n\from \invlim (Q_n,f_n)_n \to Q_n$
is a cell-like $Z^*$-map for every $n$.
\end{enumerate}
\end{theorem}

\begin{proof}
It will be convenient to let $Q_\infty$ denote $\invlim (Q_n,f_n)_n$.

By (5), every $f_n$ is a near homeomorphism. 
Hence we get (A) from Brown's Approximation Theorem for inverse limits in
\cite{brown:inverse}. 
It follows from~\cite[Theorem 6.7.4]{vm:book} that every projection
$f^{\infty}_{n}\from Q_\infty \to Q_n$ is a near homeomorphism, hence is
cell-like.

For every $n$ let $\varrho_n$ be an admissible metric for $Q_n$ which is
bounded by $1$. 
The formula
$$
\varrho(x,y)=\sum_{n=1}^\infty 2^{-n} \varrho_n(x_n,y_n)
$$
defines an admissible metric for $Q_\infty$. 
With respect to this metric we have that $f^{\infty}_{n}$ is a 
$2^{-(n-1)}$-mapping (\cite[Lemma 6.7.3]{vm:book}).

For (B) it suffices to prove that $f^{\infty}_{1}$ is a $Z^*$-map. 
To this end, let $A\con Q_1$ be a $Z$-set, and let $\varepsilon > 0$. 
Pick $n\in{\mathbb N}$ so large that $2^{-(n-1)} < \varepsilon$.
It follows that for every $x\in Q_n$ we have that the diameter of the 
fiber $(f^{\infty}_{n})^{-1}(x)$ is less than $\varepsilon$.
An easy compactness argument gives us an open cover $\mathcal{U}$ of $Q_n$ 
such that for every $U\in \mathcal{U}$ we have that
$$
\diam (f^{\infty}_{n})^{-1}[U] < \varepsilon.\eqno{(*)}
$$
Let $\gamma > 0$ be a Lebesgue number for this cover (\cite[Lemma 
1.1.1]{vm:book}).
Since $f^{\infty}_{n}$ is a near homeomorphism, there is a homeomorphism
$\varphi\from Q_\infty\to Q_n$ such that for every $x\in
Q_\infty$ we have
$$
\varrho_n(f^{\infty}_{n}(x),\varphi(x))<\leukfrac{1}{2}\gamma.
$$
Observe that $A_n=(f^{n}_{1})^{-1}[A]$ is a $Z$-set in $Q_n$. 
There consequently is a continuous function 
$\xi\from Q_n\to Q_n\setminus A_n$ which moves the points less than 
$\leukfrac{1}{2}\gamma$.
Now define $\eta\from Q_\infty \to Q_\infty$ by
$$
\eta = \varphi^{-1} \circ \xi \circ f^{\infty}_{n}.
$$
It is clear that $\eta [Q_\infty]$ misses $(f^{\infty}_{1})^{-1}[A]$.
In order to check that $\eta$ is a `small' move, pick an arbitrary element
$x\in Q_\infty$. 
By construction,
$$
\varrho_n{\big (}x_n,\xi(x_n){\big )}< \leukfrac{1}{2}\gamma.
$$
Since $\eta(x) = \varphi^{-1}{\big (} \xi(x_n) {\big )}$, clearly
$$
\varrho_n{\big (}\eta(x)_n, \xi(x_n){\big )} < \leukfrac{1}{2}\gamma.
$$ 
We conclude that $\varrho_n(\eta(x)_n,x_n)<\gamma$. 
Pick an element $U\in \mathcal{U}$ which contains both $\eta(x)_n$ and
$x_n$. 
By $(*)$ it consequently follows that $\varrho(\eta(x),x)<\varepsilon$,
which is as required.
\end{proof}

\begin{theorem}\label{tweedestelling}
If $(A_n)_n$ is a relatively discrete sequence of closed subsets of $Q$ 
such that $\CL{\bigcup_{n=1}^\infty A_n}$ is a $Z$-set then there are a
Hilbert cube $M$ and a continuous surjection $f\from M\to Q$ such that
\begin{enumerate}
\item[(A)] $f$ is a cell-like $Z^*$-map.
\item[(B)] The closures of the sets $\bigcup_{n=1}^\infty f^{-1}[A_{2n}]$
and $\bigcup_{n=0}^\infty f^{-1}[A_{2n+1}]$ are disjoint.
\end{enumerate}
\end{theorem}

\begin{proof}
Consider the subspace $A=\CL{\bigcup_{n=1}^\infty A_n}$ of $Q$, and the
`remainder' $R= A\setminus \bigcup_{n=1}^\infty A_n$.
Observe that $R$ is compact since the sequence $(A_n)_n$ is relatively
discrete.
Let $T$ denote the product $A\times \interval$; 
put
$$
S = 
(R\times \interval) 
\cup 
{\Big(}
\bigcup_{n=1}^\infty A_{2n}\times\{0\}
{\Big)} 
\cup 
{\Big(} 
\bigcup_{n=0}^\infty A_{2n+1}\times \{1\}
{\Big)}
.
$$
Then $S$ is evidently a closed subspace of $T$.
Let $\pi\from R\times\interval \to R$ denote the projection. 
It is clear that the adjunction space (cf., \cite[Page~507]{vm:book:twee})
$S\cup_\pi (R\times\interval)$ is homeomorphic to $A$. 
By (4), any constant function $S\to Q$ can be approximated by an
imbedding whose range is a $Z$-set. 
So we may assume without loss of generality that $S$ is a $Z$-subset of
some Hilbert cube $M^Q$. 
Now consider the space $N = M^Q\cup_\pi (R\times\interval)$ with natural
decomposition map $f$. 
It is clear that $f$ is cell-like, each non-degenerate fiber of $f$ being
an arc (\cite[Corollary 7.1.2]{vm:book}).
We will prove below that $N\approx Q$. 
Once we know that, we also get by Lemma~\ref{eerstelemma} that $f$ is a
$Z^*$-map.
Observe that the projection $\pi\from R\times \interval\to R$ is a 
hereditary shape equivalence. 
So by a result of Kozlowski~\cite{kozlow:anr} (see also~\cite{ancel:ce}),
it follows that $N$ is an $\mathsf{AR}$. 
Since $S$ is a $Z$-set in $M^Q$ it consequently follows 
from~\cite[Proposition 7.2.12]{vm:book} that 
$f[S]\approx A$ is a $Z$-set in $N$. 
But $N\setminus f[S]$ is obviously a $Q$-manifold, and consequently has
the disjoint-cells property. 
But this implies that $N$ has the disjoint-cells property, i.e., 
$N\approx Q$ by Toru\'nczyk's topological characterization of $Q$ 
in~\cite{tor:hilbertcube}
(see also~\cite[Corollary 7.8.4]{vm:book}). 
So we conclude that $f[S]\approx A$ is a $Z$-set in the Hilbert cube $N$. 
By (3) there is a homeomorphism of pairs $(Q,A)\approx (N,f[S])$. 
This homeomorphism may be chosen to be the `identity' on $A$. 
This shows that we are done by Lemma~\ref{eerstelemma} and the obvious 
fact that the sets
$$
\bigcup_{n=1}^\infty A_{2n}\times \{0\},\quad 
\bigcup_{n=0}^\infty A_{2n+1}\times \{1\}
$$
have disjoint closures in $M^Q$.
\end{proof}

\section{The construction}

We will now construct our example under \CH. 
After the preparatory work in \S\ref{Hilbert}, the construction is very
similar to known constructions in the literature (see e.g., 
Kunen~\cite{kunen:Lspace}).

Consider the `cube' $Q^{\omega_1}$. 
For every $1\le\alpha < \omega_1$ let $\{S_\xi^\alpha : \xi < \omega_1\}$
list all nontrivial convergent sequences in $Q^\alpha$ that do not contain
their limits.
For all $\alpha,\xi<\omega_1$ pick disjoint complementary infinite
subsets $A_\xi^\alpha$ and $B_\xi^\alpha$ of $S_\xi^\alpha$.

We shall construct for $1\le\alpha \le \omega_1$ a closed subspace
$M_\alpha\con Q^\alpha$. 
The space we are after will be $M_{\omega_1}$.

Let $\tau\from\omega_1\to\omega_1\times\omega_1$ be a surjection such that
$\tau(\beta) =\orpr{\alpha}{\xi}$ implies $\alpha\le\beta$.

For $\alpha \le \beta\le\omega_1$ let $\pi^\beta_\alpha$ be the natural 
projection from $Q^\beta$ onto $Q^\alpha$. 
The following conditions will be satisfied:
\begin{enumerate}
\item[(A)] 
$M_\alpha\approx Q$ for every $1\le \alpha<\omega_1$, and if 
$\alpha \le \beta$ then $\pi^\beta_\alpha[M_\beta] = M_\alpha$.

\noindent 
We put 
$\rho^\beta_\alpha = 
\pi^\beta_\alpha\restriction M_\beta\from M_\beta\to M_\alpha$.
\item[(B)] 
If $\alpha\le\beta$ then $\rho^\beta_\alpha\from M_\beta\to M_\alpha$ is
a cell-like $Z^*$-map.
\item[(C)] 
If $\beta<\omega_1$, $\tau(\beta)=\orpr{\alpha}{\xi}$, and
$S^\alpha_\xi\con M_\alpha$ then
$(\rho^{\beta+1}_\alpha)^{-1}[A^\alpha_\xi]$ and 
$(\rho^{\beta+1}_\alpha)^{-1}[B^\alpha_\xi]$ have disjoint closures in
$M^{\beta+1}$.
\end{enumerate}
Observe that the construction is determined at all limit ordinals 
$\gamma$. 
By compactness and (A) we must have
$$
M_\gamma = 
\{x\in Q^\gamma : 
(\forall \alpha < \gamma)(\pi^\gamma_\alpha(x)\in M_\alpha)\}.
$$
Also, if $(\gamma_n)_n$ is any strictly increasing sequence of ordinals
with $\gamma_n\nearrow\gamma$ then $M_\gamma$ is canonically homeomorphic
to
$$
\invlim {\big (}M_{\gamma_n},\rho^{\gamma_{n+1}}_{\gamma_{n}}{\big )}_n.
$$
By Theorem~\ref{eerstestelling} this implies that $M_\gamma\approx Q$ and
also that $\rho^\gamma_{\gamma_n}$ is a cell-like $Z^*$-map for every $n$. 
Since $\gamma_1$ can be \emph{any} ordinal smaller than $\gamma$, the same
argument yields that $\rho^\gamma_\alpha$ is a cell-like $Z^*$-map for
\emph{every} $\alpha < \gamma$. 
So in our construction we need only worry about successor steps.

Put $M_1 = Q^{\{0\}}$, and let $1\le\beta<\omega_1$ be arbitrary.
We shall construct $M_{\beta+1}$ assuming that $M_\beta$ has been 
constructed. 
To this end, let $\tau(\beta)=\orpr{\alpha}{\xi}$. 
We make the obvious identification of $Q^{\beta+1}$ with 
$Q^\beta\times Q$. 
If $S^\alpha_\xi\not\con M_\alpha$ then there is nothing to do. We then
fix any element $q\in Q$, and put
$$
M_{\beta+1} = M_\beta \times \{q\}.
$$
So assume that $S^\alpha_\xi\con M_\alpha$. 
By Theorem~\ref{tweedestelling} there exists a cell-like $Z^*$-map 
$f\from Q\to M_\beta$ such that
$$
f^{-1}{\big [}(\rho^\beta_\alpha)^{-1}[A^\alpha_\xi]{\big ]},
\quad
f^{-1}{\big [}(\rho^\beta_\alpha)^{-1}[B^\alpha_\xi]{\big ]}
$$
have disjoint closures in $Q$. 
Put
$$
M_{\beta+1} = \{\orpr{f(x)}{x}\in Q^{\beta}\times Q : x\in Q\}.
$$
So $M_{\beta+1}$ is nothing but the graph of $f$. 
It is clear that $M_{\beta+1}$ is as required.

Now put $M = M_{\omega_1}$. 
Observe that $M$ is a locally connected continuum, being the inverse limit
of an inverse system of locally continua with monotone surjective bonding
maps (see e.g., \cite[6.3.16 and 6.1.28]{engelking:gentop}).
Assume that $T$ is a nontrivial convergent sequence with its limit $x$ in
$M$. 
Since $T\cup \{x\}$ is countable, there exists $\alpha < \omega_1$ such
that $\rho^{\omega_1}_\beta\restriction (T\cup\{x\})$ is one-to-one and
hence a homeomorphism for every $\beta \ge \alpha$. 
Pick $\xi < \omega_1$ such that
$S^\alpha_\xi = \rho^{\omega_1}_\alpha[T]$, and $\beta \ge\alpha$ such
that $\tau(\beta) = \orpr{\alpha}{\xi}$. 
Then $\rho^{\omega_1}_{\beta+1}[T\cup\{x\}]$ is a nontrivial convergent
sequence with its limit in $M_{\beta+1}$ which is mapped by 
$\rho^{\beta+1}_\alpha$ onto $S^\alpha_\xi$ with its limit.
But this is clearly in conflict with (C).

\providecommand{\bysame}{\leavevmode\hbox to3em{\hrulefill}\thinspace}
\providecommand{\MR}{\relax\ifhmode\unskip\space\fi MR }
\providecommand{\MRhref}[2]{%
  \href{http://www.ams.org/mathscinet-getitem?mr=#1}{#2}
}
\providecommand{\href}[2]{#2}

\end{document}